\documentclass[a4paper]{amsart}

\usepackage{amsmath,amssymb}
\usepackage{amsthm,mathtools}

\newtheorem{thm}{Theorem}[section]

\newtheorem{lem}[thm]{Lemma}
\newtheorem{prop}[thm]{Proposition}

\theoremstyle{definition}
\newtheorem{defn}[thm]{Definition}
\newtheorem{exm}[thm]{Example}
\theoremstyle{remark}

\usepackage[whole]{bxcjkjatype}

\numberwithin{equation}{section}

\newcommand{\RR}{\mathbb{R}}
\newcommand{\NN}{\mathbb{N}}

\renewcommand{\d}{\delta}
\newcommand{\bp}{\mathbf{p}}

\renewcommand{\O}{\Omega}

\renewcommand{\k}{\kappa}
\newcommand{\mg}{\mathfrak{g}}

\usepackage{color}
\definecolor{darkgreen}{cmyk}{1,0,1,.2}
\definecolor{darkorchid}{rgb}{1.0, 0.5, 0}
\definecolor{persimmon}{rgb}{0.93, 0.35, 0.0}

\newdimen\theight
\def\TeXref#1{%
             \leavevmode\vadjust{\setbox0=\hbox{{\tt
                     \quad\quad  {\small \textrm #1}}}%
             \theight=\ht0
             \advance\theight by \lineskip
             \kern -\theight \vbox to
             \theight{\rightline{\rlap{\box0}}%
             \vss}%
             }}%

\begin{document}

\title[On statistics which are almost sufficient]{On statistics which are almost sufficient from the viewpoint of the Fisher metrics}

\author{Kaori Yamaguchi \qquad Hiraku Nozawa}
\address{Graduate School of Science and Engineering, Ritsumeikan University, Nojihigashi 1-1-1, Kusatsu, Shiga, 525-8577, Japan}
\email{ra0097vv@ed.ritsumei.ac.jp, hnozawa@fc.ritsumei.ac.jp}
\thanks{The first-named author is supported by JST SPRING, Grant Number JPMJSP2101. The second-named author is supported by JSPS KAKENHI Grant Numbers JP20K03620 and JP24K06723.}

\begin{abstract}
A statistic on a statistical model is sufficient if it has no information loss, namely, the Fisher metric of the induced model coincides with that of the original model due to Kullback and Ay-Jost-L\^e-Schwachh\"ofer. We introduce a quantitatively weak version of sufficient statistics such that the Fisher metric of the induced model is bi-Lipschitz equivalent to that of the original model. We characterize such statistics in terms of the conditional probability or by the existence of a certain decomposition of the density function in a way similar to characterizations of sufficient statistics due to Fisher-Neyman and Ay-Jost-L\^e-Schwachh\"ofer.
 \end{abstract}


\maketitle

\section{Introduction: Sufficient statistics and its weak version}

Sufficient statistics originally considered by Fisher \cite{Fisher} is a fundamental notion in statistics and information theory for statistical models, which are parametrized families of probability measures. According to Fisher's idea, a statistic on a statistical model is sufficient if it has all information of the model for statistical estimations on the parameters. Kullback \cite{Kullback} characterized sufficient statistics in terms of the Fisher metric, which measures the informational difference of measures in the given statistical model; a statistic is sufficient if and only if the Fisher metric of the model induced by the statistic is equal to the Fisher metric of the original model. Ay-Jost-L\^e-Schwachh\"ofer \cite{Jost, Nihat} introduced the framework of parametrized measure models to rigorously formalize the statistical models with infinite sample spaces. They redefined sufficient statistics in the context and gave its characterization in terms of density functions to generalize the result of Kullback. In this paper, we introduce almost sufficient statistics such that the Fisher metric of the induced model is bi-Lipschitz equivalent to that of the original model. Our main result is a characterization of such statistics in terms of density functions that are parallel to \cite{Nihat}. Sufficient statistics can be used to improve estimators on the statistical model by theorems of Rao-Blackwell-Kolmogorov \cite{Blackwell,Kolmogorov,Rao} and Lehmann-Scheff\'e \cite{LS}. Almost sufficient statistics can be useful for models without sufficient statistics. See \cite{YN} for the comparison of almost sufficient statistics to sufficient statistics.

Let us review quickly the framework of parametrized measure models and the definition of sufficient statistics in terms of the Fisher metric due to Ay-Jost-L\^e-Schwachh\"ofer \cite{Jost, Nihat}. Let $(M,\O,\bp)$ be a parametrized measure model, namely, $M$ is a Banach manifold, $\O$ is a measureble space and $\bp : M \to \mathcal{M}(\O)$ is a $C^1$-map, where $\mathcal{M}(\O)$ is the space of finite measures on $\O$. In this article, we always assume that $(M,\O,\bp)$ is a $2$-integrable model, for which the information loss is given by the Fisher quadratic form. Assume further that $\bp$ is dominated by a finite measure $\mu_0$ on $\O$, namely, $\bp(\xi)$ is of the form $\bp(\xi)=p(\omega;\xi)\mu_0$ for a density function $p$. Then the Fisher quadratic form on $M$ is defined by
\[
\mg(v,w) = \int_{\O} (\partial_v \log p(\cdot;\xi))(\partial_w \log p(\cdot;\xi)) d\bp(\xi)
\]
for $v,w \in T_{\xi}M$, where the right hand side is well-defined by the assumption of the $2$-integrability (see \cite{Nihat}). For a tangent vector $v$ on $M$, the \textit{(second order) information loss} of $\k$ is defined by $\mg(v,v) - \mg'(v,v)$ (\cite{Amari}), where $\mg'$ is the Fisher quadratic form of the model $(M,\Omega',\k_*\bp)$ induced by $\k$. The information loss is non-negative by the monotonicity theorem \cite{Amari,Le,Sc,H.V.} for the Fisher quadratic form;
\begin{equation}\label{eq:monot}
\mg'(v,v) \leq \mg(v,v)
\end{equation}
for every tangent vector $v$ on $M$. Due to Ay-Jost-L\^e-Schwachh\"ofer, a statistic $\k$ is defined to be sufficient if this information loss vanishes everywhere, namely, the quality holds in \eqref{eq:monot} for every tangent vector $v$ on $M$.

In this article, we consider the following quantitatively weak version of sufficient statistics. 
\begin{defn}
Let $(M,\O,\bp)$ be a $2$-integrable parametrized measure model, $\O'$ a measurable space and $0 < \delta \leq 1$. Then a statistic $\kappa : \O \to \O'$ is called $\d$-\emph{almost sufficient} for $(M,\Omega,\bp)$ if we have 
\begin{equation}\label{eq:alm}
\d^2 \mg(v,v) \leq \mg'(v,v)
\end{equation}
for all tangent vector $v$ on $M$, where $\mg$ and $\mg'$ are the Fisher quadratic forms on $M$ given by $(M,\O,\bp)$ and the model $(M,\O',\k_*\bp)$ induced by $\kappa$, respectively.
\end{defn}

By definition, a $1$-almost sufficient statistic is a sufficient statistic. By the monotonicity theorem \eqref{eq:monot}, the condition \eqref{eq:alm} is equivalent to that $\d^2 \mg(v,v) \leq \mg'(v,v) \leq \mg(v,v)$ holds for all tangent vector $v$ on $M$. If $M$ is of finite dimension, the latter implies that the Fisher-Rao distance on $M$ induced by the Fisher quadratic form $\mg'$ of the induced model is bi-Lipschitz equivalent to the original Fisher-Rao distance. 

$\d$-almost sufficient statistics can be useful for statistical models without sufficient statistics. Note that some variants of sufficient statistics, for example, Bayesian sufficient statistics \cite{BS} and linear sufficient statistics \cite{PRD,D}, have been studied in a context different from ours.

A typical example of a $\d$-almost sufficient statistic is as follows.

\begin{exm}[Binomial distributions]\label{ex:B}For $n \in \NN$, let
\begin{itemize}
    \item $M = (0,1)$,
    \item $\O = \{0,1\}^{n}$ with the counting measure and 
    \item $\bp(\xi)(\{\omega\})= \xi^{a(\omega)} (1-\xi)^{n-a(\omega)}$ for $\xi \in M$ and $\omega \in \O$, where $a(\omega)$ is the number of $1$ in the entries of $\omega$.
\end{itemize}
For $1 \leq m \leq n$, let $\Omega' = \{0, 1\}^m$ and define a statistic $\k \colon \Omega \to \Omega'$ by 
\[
\k (\omega_1, \dots, \omega_n) = (\omega_1, \dots, \omega_m).
\]
This $\k$ is $\sqrt{\frac{m}{n}}$-almost sufficient (see Section \ref{sec:ex}).
\end{exm}

Ay-Jost-L\^e-Schwachh\"ofer \cite{Nihat} characterized their sufficient statistics in terms of the conditional probability. Combining with the Fisher-Neyman characterization \cite{Neyman}, they showed that the sufficient statistics are closely related to Fisher-Neyman sufficient statistics. Let us state the characterization of sufficient statistics due to \cite{Nihat} in a detailed form.
\begin{thm}[{\cite{Kullback},\cite{Le},\cite[Chapter 5]{Nihat}, Fisher-Neyman characterization \cite{Neyman}}]\label{thm:FN0}
Let $(M,\Omega,\bp)$ be a $2$-integrable parametrized measure model of the form
\[
\bp(\xi)=p(\omega;\xi)\mu_0
\]
for a finite measure $\mu_0$ on $\Omega$ such that $p$ is positive. Let $\kappa 
:\Omega\rightarrow\Omega'$ be a statistic and let $p'(\omega;\xi)$ be the density function of $\k_*\bp(\xi)$ with respect to $\k_*\mu_0$. Assume that $M$ is connected.
Then the following are equivalent:
\renewcommand{\theenumi}{\roman{enumi}}
\renewcommand{\labelenumi}{\rm{(\theenumi)}}
\begin{enumerate}
\item $\kappa$ is sufficient for $(M,\Omega,\bp)$.
\item We have $\partial_v\log p = \partial_v\log \k^*p'$ for every tangent vector $v$ on $M$.
\item The map $\xi \mapsto \frac{p(\cdot;\xi)}{p'(\k(\cdot);\xi)}$ is a constant map.
\item \label{it:FN0-4}There exist a measurable function $s:\Omega'\times M\rightarrow\mathbb{R}$ and a function $t\in L^1(\Omega,\mu_0)$ such that we have  
\[
p(\omega;\xi)=s(\kappa(\omega);\xi)t(\omega), \quad\quad \mu_0\text{-}a.e.\ \omega \in \O \text{ and \,} \forall \xi \in M.
\]
\end{enumerate}
\end{thm}

The condition \eqref{it:FN0-4} is closely related to Fisher-Neyman sufficient statistics. Recall that a statistic $\k$ is called \emph{Fisher-Neyman sufficient} if there is a finite measure $\mu_0$ on $\O$ such that 
\[
\bp(\xi)=q(\kappa(\cdot);\xi)\mu'_0 
\] 
for some $q(\cdot;\xi)\in L^1(\Omega',\mu'_0)$, where $\mu'_0 = \k_* \mu_0$. Indeed, it is straightforward from the definition that $\k$ is Fisher-Neyman sufficient if and only if, in addition to the condition \eqref{it:FN0-4}, we have $s(\cdot;\xi)\in L^1(\Omega',\kappa_*\mu_0)$ for all $\xi \in M$ (see \cite[Theorem 5.3]{Nihat}). 

Our result is the following characterization of almost sufficient statistics parallel to Theorem \ref{thm:FN0} in a sense:
\begin{thm}\label{thm:FN}
Let $(M,\O,\bp)$ be a $2$-integrable parametrized measure model of the form
\[
\bp(\xi)=p(\omega;\xi) \mu_0
\] 
for a finite measure $\mu_0$ on $\O$ such that $p$ is positive. 
Assume that $M$ equipped with the Fisher quadratic form is a finite dimensional Riemannian manifold of class $C^2$. Let $\kappa : \O \to \O'$ be a statistic and $0< \d \leq 1$. Then the following are equivalent:
\renewcommand{\theenumi}{\roman{enumi}}
\renewcommand{\labelenumi}{\rm{(\theenumi)}}
\begin{enumerate}
\item \label{it:FN1} $\kappa$ is $\d$-almost sufficient for $(M,\O,\bp)$.
\item \label{it:FN2} We have
\begin{equation}\label{eq:FN2}
\left\|\partial_v \log \frac{p(\cdot;\xi)}{\k^*p'(\cdot;\xi)} \right\| \leq \sqrt{1-\delta^2} \| \partial_v \log p (\cdot;\xi) \|
\end{equation}
for every tangent vector $v$ on $M$, where $\|\cdot \|$ denotes the $L^2$-norm.
\item \label{it:FN3} The map $M \to L^2(\O,\mu_0); \xi \mapsto \log \frac{p(\cdot;\xi)}{p'(\k(\cdot);\xi)}$ is locally $\sqrt{1-\d^2}$-Lipschitz with respect to the Fisher-Rao distance and the distance induced by the $L^2$-metric.
\item \label{it:FN4} There exist measurable functions $s : \O'\times M \to \RR$ and $t : \O \times M \to \RR_{>0}$ such that $\log t(\cdot;\xi) \in L^2(\O,\mu_0)$ for $\forall \xi \in M$, 
\begin{itemize}
\item $p(\omega;\xi)=s(\kappa(\omega);\xi)t(\omega;\xi)$ for $\mu_0$-a.e.\ $\omega \in \O, \forall \xi \in M$  and 
\item the map $M \to L^2(\O,\mu_0); \xi \mapsto \log t(\cdot;\xi)$ is locally $\sqrt{1-\d^2}$-Lipschitz with respect to the distance defined by the Fisher quadratic form and the $L^2$-metric.
\end{itemize}
\end{enumerate}
\end{thm}

Theorem \ref{thm:FN0} holds for $k$-integrable parametrized measure model for $k$ other than $2$. However, we do not know if we can extend Theorem \ref{thm:FN} to $k$ other than $2$. Note also that, in the case where $\delta=1$, the condition \eqref{it:FN4} does not coincide with the condition (iv) of Theorem \ref{thm:FN0}.

We need to assume that $\dim M < \infty$ and the Fisher quadratic form is of class $C^2$ so that every tangent vector is tangent to a geodesic and every point of $M$ has an open neighborhood such that every pair of points can be connected by a geodesic. Note that the Fisher quadratic form of a $2$-integrable parametrized measure model is merely continuous in general. In the general situation, geodesics with given initial data may not exist (see e.g., \cite{HW}). For Riemannian manifolds of infinite dimension, it is known that there are many subtle problems even for the existence (see \cite{Ekeland}).

As the proof suggests, if every pair of points on $M$ is connected by a geodesic, we can replace the local Lipschitz condition with the Lipschitz condition.

\smallskip

\noindent \textbf{Acknowledgments.}
The authors are grateful to Masayuki Asaoka for his kind advice on examples of almost sufficient statistics. The authors also thank the members of the Saturday seminar for the stimulative discussion.

\section{Almost sufficient statistics for binomial distributions}\label{sec:ex}

Let us present some examples of almost sufficient statistics for the statistical model of binomial distributions, which is a fundamental example of an exponential family.

In the case of a statistical model with finite sample space, by definition, the Fisher metric of the induced model by a statistic $\k$ only depends on the partition of the sample space $\Omega = \sqcup_{\omega' \in \Omega'} \k^{-1}(\omega')$. Thus essentially there are only finitely different statistics for such models. It is easy to classify all sufficient statistics for exponential families by using Theorem \ref{thm:FN0}. Indeed, if $\k$ is sufficient, then for each $\omega' \in \Omega'$, all elements $\omega \in \k^{-1}(\omega')$ have a density function $p(\omega; \cdot )$ which are multiple to each other. In the case of binomial distributions, it is well known that the statistic $a \colon \Omega \to \{0,1, \dots, n\}$ that counts the number of $1$ is sufficient. Since the density functions $\xi^{i} (1-\xi)^{n-i}$ and $\xi^{j} (1-\xi)^{n-j}$ with $i \neq j$ are not multiple of each other, the statistic $a$ is minimal among sufficient statistics in the sense that any other sufficient statistic defines a subdivision of the partition given by $a$. Therefore the model has a unique minimal sufficient statistic \cite{LS}.

Let us see that there are many different examples of almost sufficient statistics. Consider the binomial distribution model $(M,\O,\bp)$, which is given by $M=(0,1)$, $\O = \{0,1\}^{n}$ with the counting measure $\mu_0$ and $\bp(\xi)(\{\omega\})= \xi^{a(\omega)} (1-\xi)^{n-a(\omega)}$ for $\xi \in M$ and $\omega \in \O$, where $a(\omega)$ is the number of $1$ in the entries of $\omega$. A typical example is Example \ref{ex:B} which is the projection $\k \colon \Omega = \{0,1\}^n \to \{0,1\}^m$ for $1 \leq m \leq n$. Indeed, since the Fisher metric $\mg$ is given by $\mg(\partial_\xi,\partial_\xi) = \frac{n}{\xi(1-\xi)}$ due to Fisher \cite{Fisher} and we can compute the Fisher metric $\mg'$ of the induced model similarly, we have  
\[
\mg'(\partial_\xi,\partial_\xi) = \| \partial_\xi \log \k^*p'( \cdot ; \xi) \|^2 = \frac{m}{\xi(1-\xi)} = \frac{m}{n} \mg(\partial_\xi,\partial_\xi),
\]
which implies that $\k$ is $\sqrt{\frac{m}{n}}$-almost sufficient, where $\partial_\xi = \frac{\partial}{\partial \xi}$ and $\| \cdot \|$ denotes the $L^2$-metric of $L^2(\O,\mu_0)$. For $m=1$, this statistic takes values in $\{0,1\}$ and hence it is minimal. In order to focus on the minimal ones, let us consider statistics with values in $\{0,1\}$, which we call \textit{a binary statistic}. Let us give a criterion of almost sufficiency of $\k$. 

\begin{prop}\label{prop:as}
  Let $\k \colon \{0,1\}^n \to \{0,1\}$ be a statistic on the model of binomial distributions (Example \ref{ex:B}). Let $f(\xi)=\k_* \bp(\xi)(\{0\})$. If $\k$ is surjective, then $\k$ is $\d$-almost sufficient for some $\d > 0$ if and only if $\partial_\xi f(\xi)$ has no zero on $[0,1]$.  
\end{prop}

We can compute the Fisher metric in terms of $f(\xi)$ as follows.
\begin{lem}
If $\kappa$ is surjective, then the Fisher metric $\mg'$ of the model induced by $\k$ is given by
\[
\mg'(\partial_\xi,\partial_\xi) = \frac{\partial_\xi f(\xi)^2}{f(\xi) \bigl(1-f(\xi) \bigr)}.
\]    
\end{lem}

\begin{proof}
Let $C_0 = \k_*\mu_0(\{0\})$ and $C_1 = \k_*\mu_0(\{1\})$, which are positive by assumption. By $f(\xi)=\k_* \bp(\xi)(\{0\}) = C_0 p'(0;\xi)$, we have $p'(0;\xi) =\frac{f(\xi)}{C_0}$. Since $\k_*\bp(\xi)$ is a probability measure on $\Omega'$, 
we have $p'(1;\xi) =\frac{1-f(\xi)}{C_1}$ and hence
\begin{multline*} 
\mg'(\partial_\xi,\partial_\xi) = \| \partial_\xi \log \kappa^* p'(\xi)\|^2 
= \\ \biggl( \partial_\xi \log \frac{f(\xi)}{C_0} \biggr)^2 f(\xi)+\biggl( \partial_\xi \log \frac{1-f(\xi)}{C_1} \biggr)^2 (1-f(\xi)) = \frac{\partial_\xi f(\xi)^2}{f(\xi) \bigl(1-f(\xi) \bigr)}.\qedhere
\end{multline*} 
\end{proof}

Consider a non-negative function on $M$
\[
h(\xi) = \frac{\mg'(\partial_\xi,\partial_\xi)}{\mg(\partial_\xi,\partial_\xi)}.
\]
By the definition of almost sufficient statistics, $\k$ is $\d$-almost sufficient for some $\d > 0$ if and only if $h(\xi)$ has a positive minimum value on $M$. We have $\mg(\partial_\xi,\partial_\xi) = \frac{n}{\xi (1-\xi)}$ due to \cite{Fisher}. Then, by the last lemma, we have 
\begin{equation}\label{eq:hxi}
h(\xi) = \frac{ \xi (1-\xi) \partial_\xi f(\xi)^2 }{nf(\xi) (1-f(\xi))}.    
\end{equation}

We will prove Proposition \ref{prop:as} with the following observation.

\begin{lem}
$\lim_{\xi \to 0} h(\xi) \neq 0$ if and only if we have $\partial_\xi f(0) \neq 0$.
\end{lem}

\begin{proof}
    Since $p(\omega; \xi)$ is a polynomial of $\xi$, so is $f(\xi)$. Take $m \geq 1$ and a polynomial $g(\xi)$ so that $f(\xi) = \xi^m g(\xi) + f(0)$ and $\xi \nmid g(\xi)$.  Then, by \eqref{eq:hxi}, we have
    \[
h(\xi) = \frac{\xi(1-\xi)(m\xi^{m-1}g(\xi) +\xi^{m}\partial_\xi g(\xi))^2}{n(\xi^m g(\xi)+f(0))(1-f(0)-\xi^{m}g(\xi))}
    \]
    For every $\omega \in \Omega$, the density function $p(\omega;\xi)$ is a polynomial and defined for $\xi \in [0,1]$. Since we have $p(\omega;0)=1$ for $\omega = (0,0,\dots,0)$ and $p(\omega;0)=0$ for $\forall \omega \in \Omega - \{(0,0,\dots,0)\}$, we have $f(0)=1$ if $\kappa((0,0,\dots,0))=0$ and $f(0)=0$ otherwise. Assume that $f(0)=1$. Then we have 
\begin{align*}
h(\xi) & = -\xi^{m-1} \frac{(1-\xi)(mg(\xi) +\xi \partial_\xi g(\xi))^2}{ng(\xi)(1+\xi^{m}g(\xi))}.
\end{align*}
Since 
\[
\frac{(1-\xi)(mg(\xi) +\xi \partial_\xi g(\xi)(\xi))^2}{g(\xi)(1+\xi^{m}g(\xi))} \to \frac{m^2 g(0)}{n} \neq 0 \quad (\xi \to 0),
\]
we have that $\lim_{\xi \to 0} h(\xi) \neq 0$ if and only if we have $m=1$, which concludes the proof of the lemma in the case where $f(0)=1$. The proof for the case where $f(0)=0$ is similar.    
\end{proof}

Similarly, we have $\lim_{\xi \to 1} h(\xi) \neq 0$ if and only if we have $\partial_\xi f(1) \neq 0$. By definition, for every $0 < \xi < 1$, we have $h(\xi)=0$ if and only if $\partial_\xi f(\xi)=0$. Therefore, Proposition \ref{prop:as} is a direct consequence of the last lemma.

Let us apply Proposition \ref{prop:as} to two examples. Since every binary almost sufficient statistic is minimal, the following example shows that there are many of them and its classification can be complicated.

\begin{exm}
Let us list up almost sufficient statistics $\k \colon \{0,1\}^n\rightarrow\{0,1\}$ for the binomial distributions for the case where $n=2,3$. In order to reduce unnecesary symmetry, we consider only $\k$ such that $\k ((1,1,\dots,1)) = 0$. In this case, $\k$ should satisfy $\k((0,0,\dots,0))=1$ to be almost sufficient. Indeed, if not, we have $\k((0,0,\dots,0))=0$, which implies $f(0) = f(1) = 1$. Then $\partial_\xi f$ should vanish at some point on $[0,1]$ by Rolle's theorem and hence $\k$ is not almost sufficient by Proposition \ref{prop:as}. 

In the case where $n=2$, for $\k$ with $\k ((0,0)) = 1$ and $\k((1,1))=0$, we have $f(\xi)=\xi^2+a \xi(1-\xi)$, where $a = \# \k^{-1}(0) \cap \{ (0,1), (1,0) \}$.  By a direct computation of $\partial_\xi f$ with Proposition \ref{prop:as}, we have that $\k$ is almost sufficient if and only if $a=1$, which essentially coincides with Example \ref{ex:B}.

In the case where $n=3$, for $\k$ with $\k ((0,0,0)) = 1$ and $\k((1,1,1))=0$, we have $f(\xi)=\xi^3+a_{1} \xi^2 (1-\xi)+a_{2} \xi(1-\xi)^2$, where $a_1 = \# \k^{-1}(0) \cap \{ (0,1,1), (1,0,1), (1,1,0) \}$ and $a_2 = \# \k^{-1}(0) \cap \{ (0,0,1), (0,1,0), (1,0,0) \}$. By a direct computation of $\partial_\xi f$ with Proposition \ref{prop:as}, we have that $\k$ is almost sufficient if and only if $(a_{1},a_{2})=(0,1),(0,2),(1,1),(1,2),(1,3),(2,1),(2,2),(2,3)$. We note that the case with $(a_1,a_2)=(2,1)$ is essentially Example \ref{ex:B} with $n=3$ and $m=1$.
\end{exm}

There are more almost sufficient statistics for binomial distrubutions for larger $n$. Since the composite of two almost sufficient statistics are almost sufficient, for the case where $n \geq 3$, we can construct binary almost sufficient statistics by the composite $\{0,1\}^n \to \{0,1\}^3 \to \{0,1\}$, where $\{0,1\}^n \to \{0,1\}^3$ is the one in Example \ref{ex:B} with $m=3$ and $\{0,1\}^3 \to \{0,1\}$ is one of the last example.

Let us mention that the majority vote is not almost sufficient.

\begin{exm}
Let $n \geq 3$ and $0\leq m \leq n-1$. Define $\kappa \colon \O \to \{0,1\}$ by
\begin{equation*}
    \kappa (\omega) = \begin{cases}
    0 & \text{ if } a(\omega) \leq m, \\
    1 & \text{ if } a(\omega) > m,
\end{cases}    
\end{equation*}
where $a(\omega)$ denotes the numbers of $1$ in the entries of $\omega$. If $m \leq n-2$, by $f(\xi)=\sum_{k=0}^{m} {}_n C_k \xi^k (1-\xi)^{n-k} =(1-\xi)^2 \sum_{k=0}^{m} {}_n C_k \xi^k (1-\xi)^{n-2-k}$, we have $\partial_\xi f(1)=0$. In the case where $m=n-1$, we have that $f(\xi) = 1 - (1-\xi)^{n}$ and $\partial_\xi f(1)=0$. Thus $\k$ is not almost sufficient by Proposition \ref{prop:as}. 
\end{exm}

\section{Characterization of $\d$-almost sufficient statistics}

Let us prove Theorem \ref{thm:FN}.
Let $(M,\O,\bp)$ be a $2$-integrable parametrized measure model with positive density function $\bp(\xi)=p(\omega;\xi) \mu_0$, where $\mu_0$ is a finite measure on $\O$.

Let us show that the conditions \eqref{it:FN1} and \eqref{it:FN2} are equivalent to each other by using an argument in \cite[Proposition 5.1]{Nihat}. Let $p'(\omega;\xi)$ be the density function of the model induced by $\k$, namely, we have $\k^*\bp(\xi) = p'(\cdot;\xi) \mu'_0$ for all $\xi$, where $\mu'_0=\k_* \mu_0$. Take a tangent vector $v$ on $M$, and let $\phi=\partial_v \log p$ and $\phi'=\partial_v \log p'$. Since $(M,\O,\bp)$ is $2$-integrable, we have $\phi \in L^2(\O,\mu_0)$ and then $\phi' \in L^2(\O',\mu'_0)$ (see, e.g., \cite[p.\ 247]{Nihat}). It is well known that the $L^2$-norm of $\phi'$ in $L^2(\O',\mu'_0)$ is equal to the $L^2$-norm of $\k^*\phi'$ in $L^2(\O,\mu_0)$ (see \cite[Eq.\ (5.14)]{Nihat}), where $\|\cdot\|$ denotes the $L^2$-norm of $L^2(\O,\mu_0)$. Therefore $\k$ is $\d$-almost sufficient if and only if $\d \| \phi \| \leq \| \k^*\phi' \|$ for all tangent vector $v$ on $M$ by definition. We have
\[
     \|\k^*\phi'\|^2
    =\|\phi'\|^2
    =\int_{\O'} (\phi')^2 d\mu'_0
    =\int_{\O'} \phi' d\k_*(\phi\mu_0)
    =\int_{\O} (\k^*\phi')\phi d\mu_0
    =\langle \phi, \k^*\phi' \rangle,
\]
which implies that
\begin{equation}\label{eq:phikphi}
\|\phi - \k^*\phi' \|^2 = \| \phi \|^2 - \| \k^*\phi' \|^2.   
\end{equation}
It follows from this equation that $\d \| \phi \| \leq \| \k^*\phi' \|$ and $\|\phi - \k^*\phi' \| \leq \sqrt{1-\d^2} \| \phi \|$ are equivalent to each other, which shows the equivalence of conditions \eqref{it:FN1} and \eqref{it:FN2}. 

Let us show the equivalence of the conditions \eqref{it:FN2} and \eqref{it:FN3}. Assume that the condition \eqref{it:FN2} holds, namely, we have 
\begin{equation}\label{eq:FN2-2}
\left\| \partial_v \log \frac{p(\cdot;\xi)}{p'(\k(\cdot);\xi)} \right\| \leq \sqrt{1-\d^2} \left\| \partial_v \log p(\cdot;\xi) \right\| 
\end{equation}
for every tangent vector $v$ on $M$. Since $\dim M < \infty$ and the Fisher quadratic form is of class $C^2$, an arbitrary point on $M$ has an open neighborhood $U$ such that every pair of points in $U$ is connected by a geodesic. Take $\forall \xi_0, \xi_1 \in U$. Let $\xi(t)$ the geodesic on $M$ with unit speed such that $\xi(t_0)=\xi_0$ and $\xi(t_1)=\xi_1$. Since the tangent vectors of $\xi(t)$ has norm one and we have $d(\xi_0,\xi_1)=|t_0-t_1|$, we have
\begin{multline*}
\left\| \log \frac{p(\cdot;\xi_0)}{p'(\kappa(\cdot);\xi_0)} - \log \frac{p(\cdot;\xi_1)}{p'(\kappa(\cdot);\xi_1)} \right\| \leq \left\| \int_{t_0}^{t_1} \frac{d}{dt} \left( \log \frac{p(\cdot;\xi(t))}{p'(\kappa(\cdot);\xi(t))} \right) dt \right\| \\
\leq \int_{t_0}^{t_1} \left\| \frac{d}{dt} \left( \log \frac{p(\cdot;\xi(t))}{p'(\kappa(\cdot);\xi(t))} \right) \right\| dt \leq \int_{t_0}^{t_1} \sqrt{1-\d^2} dt = \sqrt{1-\d^2} d(\xi_0,\xi_1).
\end{multline*}
This means that $M \to L^2(\O,\mu_0); \xi \mapsto \log \frac{p(\cdot;\xi)}{p'(\kappa(\cdot);\xi)}$ is $\sqrt{1-\d^2}$-Lipschitz on $U$, which shows that the condition \eqref{it:FN3} holds. 

Assume conversely that $\log \frac{p(\cdot;\xi)}{p'(\kappa(\cdot);\xi)}$ is locally $\sqrt{1-\d^2}$-Lipschitz.
Take an arbitrary $\xi_0 \in M$ and $v \in T_{\xi_0}M$. It is sufficient to consider the case where the norm of $v$ is one, which we will consider below. Since the Fisher quadratic form is of class $C^2$, there exists a geodesic $\xi(t)$ on $M$ such that $\xi (0)=\xi_0$ and $\frac{d\xi}{dt}(0)=v$. By assumption, for $t$ sufficiently close to $0$, we have
\[
\frac{1}{d(\xi_0,\xi(t))}\left\| \log \frac{p(\cdot;\xi)}{p'(\kappa(\cdot);\xi)} - \log \frac{p(\cdot;\xi(t))}{p'(\kappa(\cdot);\xi(t))} \right\| \leq \sqrt{1-\d^2} 
\]
By $d(\xi_0,\xi(t))=t$, we get \eqref{eq:FN2-2} by taking the limit $t \to 0$.

Finally let us show the equivalence of the conditions \eqref{it:FN3} and \eqref{it:FN4}. If the condition \eqref{it:FN3} holds, then the condition \eqref{it:FN4} is satisfied with $s=\kappa^*p'$ and $t=\frac{p}{\k^*p'}$. Conversely, we assume that the condition \eqref{it:FN4} holds, namely, there exist measurable functions $s : \O'\times M \to \RR$ and $t : \O \times M \to \RR_{>0}$ such that $p(\omega;\xi)=s(\k(\omega);\xi)t(\omega;\xi)$ for $\mu_0$-a.e.\ $\omega \in \O, \forall \xi \in M$ and the map $M \to L^2(\Omega,\mu_0); \xi \mapsto \log t(\cdot;\xi)$ is locally $\sqrt{1-\d^2}$-Lipschitz.

Let $t'(\cdot;\xi)$ be the measurable function on $\O'$ such that $\k_*(t(\cdot;\xi)\mu_0)=t'(\cdot;\xi)\mu'_0$. Note that, since $p$ is positive by assumption, $s(\cdot,\xi)$ and $p'(\k(\cdot),\xi)$ are positive as well. By
\[
p'(\kappa(\cdot);\xi)\mu'_0=\bp'(\xi)=\kappa_*\bp(\xi)=
s(\k(\cdot);\xi)t'(\k(\cdot);\xi))\mu'_0,
\]
we have $t'(\k(\cdot);\xi_0)= \frac{s(\k(\cdot);\xi_1)}{p'(\k(\cdot);\xi_1)}$. By assumption, for every point on $M$, there exists an open neighborhood $U$ such that the map $U \to L^2(\O,\mu_0); \xi \mapsto \log t(\cdot;\xi)$
  is $\sqrt{1-\d^2}$-Lipschitz. Then, for all $\xi_0,\xi_1 \in U$, we have
\begin{align*}
    & \left\| \log \frac{p(\cdot;\xi_0)}{p'(\kappa(\cdot));\xi_0)}-\log\frac{p(\cdot;\xi_1)}{p'(\kappa(\cdot);\xi_1)}\right\|^2 \\
    = & \left\| \log \left(\frac{p(\cdot;\xi_0)}{s(\kappa(\cdot);\xi_0)}\cdot\frac{s(\kappa(\cdot);\xi_0)}{p'(\kappa(\cdot);\xi_0)}\right)-\log\left(\frac{p(\cdot;\xi_1)}{s(\kappa(\cdot);\xi_1)}\cdot\frac{s(\kappa(\cdot);\xi_1)}{p'(\kappa(\cdot);\xi_1)}\right)\right\|^2 \\
    = & \left\| \log t(\cdot;\xi_0) - \log t(\cdot;\xi_1)-\k^*\big(\log t'(\cdot;\xi_0)- \log t'(\cdot;\xi_1)\big)\right\|^2\\
    = & \| \log t(\cdot;\xi_0) -\log t(\cdot;\xi_1)\|^2-\| \k^*(\log t'(\cdot;\xi_0)- \log t'(\cdot;\xi_1))\|^2\\
    \leq & \|\log t(\cdot;\xi_0)-\log t(\cdot;\xi_1)\|^2 \leq (1-\delta^2)d(\xi_0,\xi_1)^2,
\end{align*}
where the third equality follows in the same way as \eqref{eq:phikphi}. This means that the map $\xi \mapsto \log \frac{p(\cdot;\xi)}{p'(\kappa(\cdot));\xi)}$ is $\sqrt{1-\d^2}$-Lipschitz on $U$, namely, the condition \eqref{it:FN3} holds.

\end{document}